\documentclass{article}
\usepackage{amssymb}
\usepackage{amsfonts}
\usepackage{amsmath}
\usepackage{hyperref}
\usepackage{graphicx}
\hypersetup{urlcolor=blue,linkcolor=blue ,citecolor=blue,colorlinks=true}
\newtheorem{theorem}{Theorem}[section]

\newtheorem{definition}[theorem]{Definition}
\newtheorem{example}[theorem]{Example}

\newtheorem{lemma}[theorem]{Lemma}

\newtheorem{remark}[theorem]{Remark}

\newenvironment{proof}[1][Proof]{\textbf{#1.} }{\ \rule{0.5em}{0.5em}}

\begin{document}

\title{Generalization  of the final-value theorem and its application in fractional differential systems}
\author{ Yayun Wu\\
{\small  \textit{Beijing Normal University, Beijing 100875,
China} }}
\maketitle

\section{Abstract}
In this paper, we generalized the known Laplace-transform final-value theorem. From our conclusion, one can deduce the existing results in $[1,3,12]$. By using final value theorem, we give a new proof that Caputo fractional differential equations have no nonconstant periodic solution.
\section{Introduction}
The final value theorem is an extremely useful result in Laplace transform theory. The final value theorem provides an explicit technique for determining the asymptotic value of a signal without having to first invert the Laplace transform to determine the time signal$[1]$.The theorem have been applied in control systems,queuing theory,ergodic physical$[6]$.

The standard assumptions for the final value theorem$[2]$ require that the Laplace transform result have all of its poles either in the open left half plane or at the origin,with at most a single pole. In this condition, the time  function has a finite limit, then one can have
\begin{equation*}
\lim_{s\rightarrow 0}sF(s)=\lim_{t\rightarrow\infty}f(t),\eqno(1)
\end{equation*}
where $F(s)=\mathcal{L}\{f\}(s)$(see definition $3.1$) and  $s$ approach zero through the right half plane.

In the standard final value theorem ,we need both side limit in $(1)$ should exist. In $[1]$, the authors publicize and prove the "infinite-limit" version of the final value theorem which can be applied to irrational functions well.
In $[3]$,the authors consider a function $f(t)$ which is periodic or asymptotically equal to a sum of periodic functions, then
\begin{equation*}
\lim_{s\rightarrow 0}sF(s)=\left \langle f\right\rangle=\lim_{t\rightarrow\infty}\frac{\int_{0}^{t}f(u)du}{t}.\eqno(2)
\end{equation*} Clearly,$(2)$ is a generalization of $(1)$. But for some function like $f(t)=t^{2}\sin(t)$(see example $4.3$), $(1)$ and $(2)$ both  can't work.

In this paper,we general the results in $[3]$ to the general case which can be used to more functions. The existence of periodic solution is a desired property in fractional dynamical systems. By using the final value theorem of periodic function, we can prove that Caputo fractional differential equations  have no nonconstant periodic solution.
\section{Preliminaries}
In this section,we introduce some definitions and results which will be useful throughout the paper. One can see $[4,5,13]$ for more details. In this paper, we always assume that $s\in \mathbb{C}$ and $s\rightarrow0$ means $s$ approach zero through the right half plane.
\begin{definition}If $f$ is piecewise continuous on $[0,\infty)$ and of exponential order $c>0$,then Laplace transform $\mathcal{L}\{f\}(s)$ exists for $\Re(s)>c$ and converges absolutely.
\begin{equation*}
F(s)=\mathcal{L}\{f\}(s)=\int_{0}^{\infty}e^{-st}f(t)dt.\eqno(3)
\end{equation*}
\end{definition}

\begin{definition}
Let $\alpha>0,[a,b]\subset \mathbb{R}$ and $x:[a,b]\rightarrow \mathbb{R}$ be a measurable function and such that $\int_{a}^{b}|x(\tau)|d\tau<\infty$.The $Riemman-Liuville$ fractional integral operator of order $\alpha$  $$(_{a}I_{t}^{\alpha}x)(t)=\frac{1}{\Gamma(\alpha)}\int_{a}^{t}(t-\tau)^{\alpha-1}x(\tau)d\tau. \eqno(4)
$$
in which $\Gamma(\cdot)$ is Gamma function. Set $_{a}I_{t}^{0}x(t)=x(t)$. From the definition,one can have
\begin{equation*}
{_{a}I_{t}^{\alpha}}\cdot{_{a}I_{t}^{\beta}}f(t)={_{a}I_{t}^{\alpha+\beta}}f(t).\eqno(5)
\end{equation*}
\end{definition}
\begin{definition}
The $\alpha$-order Caputo's-fractional derivative is defined by
 $$(^{c}_{a}D_{t}^{\alpha}x)(t)=
{_{a}I}_{t}^{m-\alpha}\frac{d^{m}x(t)}{dt^{m}},
\eqno(6)
$$
where $m=\lceil\alpha\rceil,m\in\mathbb{N}^{+}.$ $x(t)\in L^{1}([a,t]).$ When $\alpha=m,m\in\mathbb{N},$ define $D_{t}^{m}=\frac{d^{m}}{dt^{m}}$.
\end{definition}
\begin{lemma}{$([5,p.61])$}
Let $f(t)$ is continuous and exponential bounded in $[0,\infty)$, $F(s)=\mathcal{L}\{f\}(s)$  and  $s_{0}$ is a point in the region of convergence of $F(s)$, $l$ is a positive number,and if
\begin{equation*}
F(s_{0}+nl)=\int_{0}^{\infty}e^{-(s_{0}+nl)t}f(t)dt=0,(n=0,1,2,\cdots),\eqno(9)
\end{equation*}
then
\begin{equation*}
f(t)\equiv 0 .
\end{equation*}
\end{lemma}
\section{Main Result}
In this part,we will prove the generalized Laplace final value theorem.  We give an example to illustrate  our results.

\begin{theorem} Let $f(t) \in L^{1}_{loc}(\mathbb{R}^{+})$ and $F(s)=\mathcal{L}\{f\}(s)$. If $F(s)$ doesn't have poles in the open right half plane and $\lim\limits_{s\rightarrow0}sF(s)$ exists, then $\forall \alpha\geq0$
\begin{equation*}
\lim_{s\rightarrow 0}\frac{sF(s)}{\alpha+1}=\lim_{s\rightarrow 0}s\mathcal{L}\{\frac{\Gamma(\alpha+1){_{0}I_{t}^{\alpha+1}f(t)}}{t^{\alpha+1}}\}(s),\eqno(10)
\end{equation*}
and if $\exists \alpha_{0}\geq0,$ s.t. $\lim\limits_{t\rightarrow\infty}\frac{\Gamma(\alpha_{0}+1){_{0}I_{t}^{\alpha_{0}+1}f(t)}}{t^{\alpha_{0}+1}}$ exists,then
\begin{equation*}
\lim_{s\rightarrow 0}\frac{sF(s)}{\alpha_{0}+1}=\lim_{t\rightarrow \infty}\frac{\Gamma(\alpha_{0}+1){_{0}I_{t}^{\alpha_{0}+1}f(t)}}{t^{\alpha_{0}+1}}.\eqno(11)
\end{equation*}

\begin{proof}$(I)$ In fact, we have
\begin{equation*}
\mathcal{L}\{\frac{\Gamma(\alpha+1){_{0}I_{t}^{\alpha+1}f(t)}}{t^{\alpha+1}}\}(s)=\int_{1}^{\infty}\frac{F(su)}{u}(1-\frac{1}{u})^{\alpha}du,
\end{equation*}
then we only need to prove
\begin{equation*}
\lim_{s\rightarrow0}s\int_{1}^{\infty}\frac{F(su)}{u}(1-\frac{1}{u})^{\alpha}du=\lim_{s\rightarrow0}\frac{sF(s)}{\alpha+1}.
\end{equation*}

$(i)$ If $\lim\limits_{s\rightarrow0}F(s)=\infty$,and $\lim\limits_{s\rightarrow0}sF(s)$ exists,we can know that $s=0$ is the simple pole of $F(s)$, for $\forall \gamma\geq0$
\begin{equation*}
\lim_{s\rightarrow 0}\int_{1}^{\infty}\frac{F(su)}{u^{\gamma+1}}du=\infty,\lim_{s\rightarrow 0}\int_{1}^{\infty}\frac{|sF(su)|}{u}du<\infty.\eqno(12)
\end{equation*}

$(a)$ when $\alpha$ is non-negative integer number,
\begin{equation*}
\int_{1}^{\infty}\frac{F(su)}{u}(1-\frac{1}{u})^{\alpha}du=\sum_{j=0}^{\alpha}C_{\alpha}^{j}(-1)^{j}\int_{1}^{\infty}\frac{F(su)}{u^{j+1}}du
\end{equation*}
hence
\begin{equation*}
\lim_{s\rightarrow 0}s\int_{1}^{\infty}\frac{F(su)}{u}(1-\frac{1}{u})^{\alpha}du=\sum_{j=0}^{\alpha}C_{\alpha}^{j}\frac{(-1)^{j}}{j+1}\lim_{s\rightarrow0}sF(s)=\lim_{s\rightarrow0}\frac{sF(s)}{\alpha+1},
\end{equation*}
in which $\lim\limits_{s\rightarrow0}s\int_{1}^{\infty}\frac{F(su)}{u^{j+1}}du=\lim\limits_{s\rightarrow0}s^{j+1}\int_{s}^{\infty}\frac{F(v)}{v^{j+1}}dv=\lim\limits_{s\rightarrow0}\frac{sF(s)}{j+1}$ is used.

$(b)$ when $\alpha$ is a positive fraction,
\begin{equation*}
\int_{1}^{\infty}\frac{F(su)}{u}(1-\frac{1}{u})^{\alpha}du=\int_{1}^{\infty}\frac{F(su)}{u}[1+\sum_{n=1}^{\infty}\frac{(-\alpha)_{n}}{n!}\frac{1}{u^{n}}]du,
\end{equation*}
where $(-\alpha)_{n}=\prod\limits_{j=0}^{n-1}(j-\alpha)=\frac{\Gamma(n-\alpha)}{\Gamma(-\alpha)}$.
Since
\begin{equation*}
|\Gamma(z+a)/\Gamma(z+b)|\sim |z|^{a-b}, |z|\rightarrow\infty,|\arg(z)|<\pi-\varepsilon,\varepsilon>0.\eqno(13)
\end{equation*}
$([7,p.33])$ and $u>1$, one can find  $K>0$, so that
\begin{equation*}
\sum_{n=1}^{\infty}|\frac{(-\alpha)_{n}}{n!}\frac{1}{u^{n}}|\leq \sum_{n=1}^{\infty}\frac{K}{n^{\alpha+1}}.
\end{equation*}
Since $\sum\limits_{n=1}^{\infty}\frac{K}{n^{\alpha+1}}$ exists and with $(12)$, we can exchange the order of calculus and limit, then
\begin{equation*}
\lim_{s\rightarrow0}s\int_{1}^{\infty}\frac{F(su)}{u}(1-\frac{1}{u})^{\alpha}du
=\lim_{s\rightarrow0}sF(s)\sum_{n=0}^{\infty}\frac{(-\alpha)_{n}}{n!}\frac{1}{n+1}=\lim_{s\rightarrow0}\frac{sF(s)}{\alpha+1}.
\end{equation*}

$(ii)$ If $\lim\limits_{s\rightarrow 0}|F(s)|<\infty,$ then $\lim\limits_{s\rightarrow 0}sF(s)=0.$ According with $\lim\limits_{\Re s\rightarrow\infty}F(s)=0$,we have
\begin{equation*}
\lim_{s\rightarrow0}\int_{1}^{\infty}\frac{F(su)}{u}du=\infty,\lim_{s\rightarrow0}\int_{1}^{\infty}\frac{|F(su)|}{u^{\gamma}}du<\infty,\gamma>1.
\end{equation*}
By using the same method with $(i)$,we can acquire $(10)$.

$(II)$ When $\lim\limits_{t\rightarrow\infty}\frac{\Gamma(\alpha_{0}+1){_{0}I_{t}^{\alpha_{0}+1}f(t)}}{t^{\alpha_{0}+1}}$ exists, then by using $(1)$,we have
\begin{equation*}
\lim_{s\rightarrow 0}\frac{sF(s)}{\alpha_{0}+1}=\lim_{t\rightarrow \infty}\frac{\Gamma(\alpha_{0}+1){_{0}I_{t}^{\alpha_{0}+1}f(t)}}{t^{\alpha_{0}+1}}.
\end{equation*}
\end{proof}
\end{theorem}

\begin{remark}
In many applications, it may be hard to find the  $\alpha_{0}$.But if $f(t)$ is a bounded, periodic function and of period $T>0$, we can choose $\alpha_{0}=0$ as the smallest number which can make $(11)$ success. In fact,we denote $g(t)=f(t)-\frac{1}{T}\int_{0}^{T}f(\tau)d\tau$,one can easily check that $ \int_{a}^{a+T}g(\tau)d\tau=0,\forall a\geq0.$  When $t\in((n-1)T,nT],n=\mathbb{N}^{+},$
\begin{equation*}
|\frac{\int_{0}^{t}g(\tau)d\tau}{t}|=|\frac{\int_{(n-1)T}^{t}g(\tau)d\tau}{t}|\leq\frac{\int_{0}^{T}|g(\tau)|d\tau}{t}\rightarrow0, t\rightarrow\infty.
\end{equation*}
which means
\begin{equation*}
\lim_{t\rightarrow\infty}\frac{\int_{0}^{t}f(\tau)d\tau}{t}=\frac{\int_{0}^{T}f(\tau)d\tau}{T}.
\end{equation*}
Hence
\begin{equation*}
\lim_{s\rightarrow0}sF(s)=\frac{\int_{0}^{T}f(\tau)d\tau}{T}.\eqno(14)
\end{equation*}
If $f(t)$ is also differentiable over a period,then
\begin{equation*}
\lim_{s\rightarrow 0}\mathcal{L}\{f'(t)\}(s)=\frac{\int_{0}^{T}f(t)dt}{T}-f(0^{+}).\eqno(15)
\end{equation*}
\end{remark}

\begin{example} Let $f(t)=t^{q}\sin(\omega t),q\geq0,\omega>0$, $F(s)=\mathcal{L}\{f\}(s)$, then for $\forall\alpha\geq q$, we have
\begin{equation*}
\lim_{s\rightarrow0}\frac{sF(s)}{\alpha+1}=\lim_{t\rightarrow\infty}\frac{\Gamma(\alpha+1){_{0}I_{t}^{\alpha+1}}f(t)}{t^{\alpha+1}}=0.
\end{equation*}
\begin{proof}Let $f(t)=t^{q}\sin(\omega t)$, then by $(3)$,one can have
\begin{equation*}
F(s)=\frac{\Gamma(q+1)[(s+i\omega)^{q+1}-(s-i\omega)^{q+1}]}{2i(s^{2}+\omega^{2})^{q+1}}.
\end{equation*}
$F(s)$ doesn't have poles in the open right half plane
 and $\lim\limits_{s\rightarrow0}sF(s)=0.$ By using theorem $4.3$, there exists $\alpha\geq 0$, so that
\begin{equation*}
\lim_{t\rightarrow\infty}\frac{\Gamma(\alpha+1)}{t^{\alpha+1}}{_{0}I_{t}^{\alpha+1}}t^{q}\sin(\omega t)=0.
\end{equation*}
In fact,we  choose $\alpha\geq q$. By using the formula $3.385$ in $[9,p.349]$,one can have
\begin{equation*}
\frac{\Gamma(\alpha+1)}{t^{\alpha+1}}{_{0}I_{t}^{\alpha+1}}t^{q}\sin(\omega t)=\int_{0}^{1}(1-u)^{\alpha}(tu)^{q}\sin(\omega tu)du
\end{equation*}
\begin{equation*}
=B(\alpha+1,q+2)\omega t^{q+1}\sum_{n=0}^{\infty}\frac{(\frac{q}{2}+1)_{n}(\frac{q+3}{2})_{n}}{(\frac{3}{2})_{n}(\frac{q+\alpha+3}{2})_{n}(\frac{q+\alpha}{2}+2)_{n}}\frac{(-\frac{\omega^{2}t^{2}}{4})^{n}}{n!},
\end{equation*}
in which $B(\cdot,\cdot)$ is beta function$[7]$.
According with $(13)$,
\begin{equation*}
\frac{(\frac{q}{2}+1)_{n}(\frac{q+3}{2})_{n}}{(\frac{3}{2})_{n}(\frac{q+\alpha+3}{2})_{n}(\frac{q+\alpha}{2}+2)_{n}}\leq\frac{1}{n^{\alpha+1}}\leq1,n\rightarrow\infty.
\end{equation*}
so one can  find proper $M>0$, so that
\begin{equation*}
\sum_{n=0}^{\infty}\frac{(\frac{q}{2}+1)_{n}(\frac{q+3}{2})_{n}}{(\frac{3}{2})_{n}(\frac{q+\alpha+3}{2})_{n}(\frac{q+\alpha}{2}+2)_{n}}\frac{(-\frac{\omega^{2}t^{2}}{4})^{n}}{n!}\leq Me^{-\frac{\omega^{2}t^{2}}{4}},
\end{equation*}
hence
\begin{equation*}
\lim_{t\rightarrow\infty}|\frac{\Gamma(\alpha+1)}{t^{\alpha+1}}{_{0}I_{t}^{\alpha+1}}t^{q}\sin(\omega t)|\leq\lim_{t\rightarrow\infty}\frac{MB(\alpha+1,q+2)\omega t^{q+1}}{e^{\frac{\omega^{2}t^{2}}{4}}}=0.
\end{equation*}
\end{proof}
\end{example}

\section{Nonexistence of nonconstant periodic solutions in fractional-order dynamical systems}

The existence of periodic solutions is a desired property in fractional dynamical systems. Some results have been investigated in $[10,11]$. In this paper, we give a new method to study the period solutions of the following scalar fractional system
\begin{equation*}
\begin{cases}
{^{c}_{0}D_{t}^{\alpha}x(t)=f(x(t))},t\geq0,\alpha\in(0,1).\\
x(0)=x_{0}.
\end{cases}\eqno(16)
\end{equation*}

\begin{theorem}
 Fractional system $(16)$ has no nonconstant periodic solution.

\begin{proof}
If $x(t)$ is a nonconstant periodic solution of $(16)$ and of period $T>0$,which means
$^{c}_{0}D^{\alpha}_{t}x(t)$ will be periodic function with same period, then
\begin{equation*}
\mathcal{L}\{^{c}_{0}D^{\alpha}_{t}x(t)\}(s)=\frac{1}{s^{1-\alpha}}\mathcal{L}\{x'(t)\}(s).
\end{equation*}
By using $(14)$ and $(15)$,
\begin{equation*}
\Rightarrow 0=\lim_{s\rightarrow 0}s^{\alpha}\mathcal{L}\{x'(t)\}(s)=\lim_{s\rightarrow 0}s\mathcal{L}\{^{c}_{0}D^{\alpha}_{t}x(t)\}(s)=\frac{\int_{0}^{T}{^{c}_{0}D^{\alpha}_{t}x(t)}dt}{T}
\end{equation*}
which means
\begin{equation*}
\int_{0}^{T}(T-\tau)^{1-\alpha}x'(\tau)d\tau=0.
\end{equation*}
let $\tau=T(1-e^{-\frac{(n+1)u}{1-\alpha}}),n=0,1,2,\cdots$, then
\begin{equation*}
\int_{0}^{\infty}e^{-(n+1)\frac{2-\alpha}{1-\alpha}u}x'(T(1-e^{-\frac{(n+1)u}{1-\alpha}}))du=0.
\end{equation*}
If we denote $s_{0}=l=\frac{2-\alpha}{1-\alpha}$, by using lemma $3.4$,one can know that $x'(T(1-e^{-\frac{(n+1)t}{1-\alpha}}))\equiv 0 $,which means $x(t) \equiv$ constant. That would be a contradiction with $x(t)$ is a nonconstant periodic solution. So fractional system $(16)$ has no nonconstant periodic solution.
\end{proof}
\end{theorem}
 \begin{remark}
 From the definition of Caputo derivative, one can easily check that when $n-1<\alpha<n,n\in\mathbb{Z}^{+}$, the theorem $5.1$ can also be true. By using this method, the conclusion is also true for the multidimensional fractional systems.
\end{remark}

\end{document}